\documentclass[10pt,oneside,a4paper]{article}
\usepackage{amsmath,amsthm}
\usepackage{hyperref}
\usepackage{geometry}
\usepackage{url}
\usepackage{amssymb}

\newtheoremstyle{mytheorem}
  {3pt}
  {3pt}
  {\itshape}
  {}
  {\bfseries}
  {.}
  {1em}
  {}
\theoremstyle{mytheorem}
\newtheorem*{theorem}{Theorem}
\newtheorem*{thm}{Theorem (Liu)}

\newtheorem*{definition}{Definition}
\newtheorem*{def 1}{Definition (Edwards)}
\newtheorem*{def 2}{Definition ($\varepsilon$-distance)}

\newtheorem*{proposition}{Proposition}

\usepackage{sectsty}
\sectionfont{\centering}
\usepackage{secdot}
\newtheorem{lemma}{Lemma}

\newtheorem*{corollary}{Corollary}

\theoremstyle{remark}
\newtheorem{remark}{Remark}

\begin{document}
\title{Metric topology on the moduli space}
\author{Jialong Deng}
\date{}
\newcommand{\Addresses}{{
  \bigskip
  \footnotesize
  \textsc{Mathematisches Institut, Georg-August-Universität, Göttingen, Germany}\par\nopagebreak
  \textit{E-mail address}: \texttt{jialong.deng@mathematik.uni-goettingen.de}}}
\maketitle
\begin{abstract}
We define the smooth Lipschitz topology on the moduli space and show that each conformal class is dense in the moduli space endowed with Gromov-Hausdorff topology, which offers an answer to the Tuschmann's question. 
\end{abstract}\par

 Any smooth closed manifold $M$ can be given a smooth Riemannian metric, and then one can ask: How many Riemannian metrics are there, and how many different geometries of this kind does the manifold actually allow? That means  one want to understand the space of Riemannian metrics on $M$, which is denoted by $\mathcal{R}(M)$, and the moduli space which is denoted by $\mathcal{M}(M)$. Here the moduli space is  the quotient space of $\mathcal{R}(M)$ by the action of diffeomorphism group of $M$. Those two questions originated from Riemann when he set up the Riemannian geometry in the nineteenth century \cite{2016arXiv160207208A}. Especially, the moduli space $\mathcal{M}(M)$ is the superspace in the physics [see \cite{MR0343855}, \cite{MR0347323}, \cite{MR0401069}].\par
   The $C^{n,\alpha}$-compact-open topology ($n\in \mathbb{N}^+$ and $\alpha \in \mathbb{R}^+$) is the most common consideration \cite{MR3445334}. Tuschmann asked the following question in [\cite{MR3547932}, Section 3, (8)]: What can one say about the topology of moduli spaces under the Gromov-Hausdorff metric? What if one uses the Lipschitz topology?  Inspired by Tuschmann's questions, we will introduce four kinds of metric topology on the moduli space, and then discuss the relationship among them. \par

Let $X$ and  $Y$ are metric spaces of finite diameter, then the Gromov-Hausdorff distance is defined by $\rho_{GH}(X,Y):=\inf \limits_{Z}\{d_H^Z(f(X), g(Y))\} $ where $d_H$ is Hausdorff metric and $Z$ takes all metric spaces such that $f$ (resp. $g$)  are isometric embeddings $X$ (resp. $Y$) into $Z$ [See, \cite{MR623534}]. The  Gromov-Hausdorff distance $\rho_{GH}$ is a pseudo-metric on the collection of all compact metric spaces. Furthermore, $\rho_{GH}(X, Y)=0$ if and only if $X$ is isometric to $Y$.  For $g_1$ and $g_2$ in $\mathcal{R}(M)$, the Gromov-Hausdorff distance can be defined  on it by  $\rho_{GH}(g_1, g_2)=\rho_{GH}((M,d_{1}),(M,d_{2}))$ where $d_1$ and $d_2$ are induced metrics on $M$ by $g_1$ and $g_2$. Since $M$ is closed, the Gromov-Hausdorff  distance is well-defined on $\mathcal{R}(M)$. Moreover, $\rho_{GH}(f_1^*g_1, f_2^*g_2)=\rho_{GH}(g_1, g_2)$ where $f_1$ and $f_2$ are  diffeomorphism of $M$ and $f_1^*g_1$, $f_2^*g_2$ are push-back metrics on $M$. Then one can define $\rho_{GH}$  on $\mathcal{M}(M)$ as above and then $\rho_{GH}$ is a metric on $\mathcal{M}(M)$. Therefore, $\mathcal{M}(M)$ can be endowed with the metric topology  called GH-topology by the Gromov-Hausdorff metric  $\rho_{GH}$.

\begin{def 1}
 A map $f: X\rightarrow Y$ is called an $\varepsilon$-isometry between compact metric spaces $X$ and $Y$, if   $\left|d_{X}(a, b)-d_{Y}(f(a), f(b))\right|\leq\varepsilon$ for all $a,b \in X$.
\end{def 1}
 
 \begin{def 2}
  Let $g_1$ and $g_2$ in $ \mathcal{R}(M)$, then the $\varepsilon$-distance is defined by   $\rho_\varepsilon(g_1, g_2):=\rho_\varepsilon (d_1, d_2)=\inf\left\{\varepsilon\mid I_{\varepsilon}(d_1, d_2)\neq\phi\neq I_{\varepsilon}(d_2, d_1)\right\} $, where   $I_{\varepsilon}(d_1, d_2)$ be the set of $\varepsilon$-isometries from $(M, d_1)$ to $(M, d_2)$.
\end{def 2}
  
\begin{remark}
   Since $\left|\mathrm{Diam}(d_1)-\mathrm{Diam}(d_2)\right|\leq\rho_{\varepsilon}(d_1, d_2)\leq  \max\{\mathrm{Diam}(d_1),\mathrm{Diam}(d_2)\}$,  where $\mathrm{Diam}(d_i)$ is the diameter of $(M,g_i)$, $i=1,2$,  $\rho_{\varepsilon}$ is well-defined on $\mathcal{R}(M)$. Moreover, $\rho_{\varepsilon}$ is the pseudo-metric  and $\rho_{\varepsilon}(g_{1}, g_{2})=0$ if and only if $g_{1}$ is isometric to $g_{2}$ on $\mathcal{R}(M)$. 
\end{remark}

 Then $\varepsilon$-metric, which also be denoted by $\rho_{\varepsilon}$,  can be defined on $\mathcal{M}(M)$ as $\rho_{GH}$. Thus, it induces a metric topology  on $\mathcal{M}(M)$  called $\varepsilon$-topology. The conformal class dense theorem of the $\varepsilon$-topology on $\mathcal{M}(M)$ was proved by  Liu in [\cite{2015arXiv150903895L}, Corollary 2.2].
\begin{thm} 
Each conformal class is dense in $\mathcal{M}(M)$ which is endowed with $\varepsilon$-topology.
\end{thm}

\begin{lemma}
 If $\rho_{GH}(X,Y)\leq\varepsilon$, then there is a 2$\varepsilon$-isometric map $f:X\rightarrow Y$. If there is an $\varepsilon$-isometric map $f:X\rightarrow Y$, then $\rho_{GH}(X,Y)\leq\frac{3}{2}\varepsilon$. 
\end{lemma}

\begin{remark}
The lemma can be proved by using another definition of Gromov-Hausdorff metric, i.e. $ \rho_{GH}(X,Y)=\frac{1}{2}\underset{R}{\mathrm{inf}}\left\{ \mathrm{dis}\left(R\right)\right\}$, where the infimum is taken over all correspondences $R\subseteq X\times Y$. A correspondence  between two metric spaces $X$ and $Y $ is a subset $R$ of $X\times Y$ such that the projections $\pi_X: X\times Y \rightarrow X$ and $\pi_Y: X\times Y\rightarrow Y$ remain surjective when they are restricted to $R$. 
\end{remark}

 \begin{corollary}
 $\varepsilon$-topology is  equivalent to GH-topology on $\mathcal{M}(M)$.
\end{corollary}
Thus, the conformal class dense theorem is also ture on GH-topology.  That means the GH-topology is coarse and the more finer topology is needed to define on $\mathcal{M}(M)$.\par

Let $X$ and $Y$ are two compact metric spaces,  the dilation of a Lipschitz map $f:X\rightarrow Y$ is defined by $\mathrm{Dil}(f):=\sup\limits_{a,b\in X, a\neq b} \frac{d_Y(f(a),f(b))}{d_X(a,b)}$. If $f^{-1}$ is also a Lipschitz map then it is called the bi-Lipschitz homeomorphism. The Lipschitz-distance $\rho_L$ between $X$ and $Y$ is defined by $\rho_L(X,Y):= \inf \limits_{f:X\rightarrow Y} \log\{\max\{\mathrm{Dil}(f),\mathrm{Dil}(f^{-1})\}\}$ where the infinum is taken over all bi-Lipschitz homeomorphisms between $X$ and $Y$.  Then  the Lipschitz-distance $\rho_L$ can be defined on $\mathcal{R}(M)$ as the definition of Gromov-Hausdorff distance on $\mathcal{R}(M)$.\par
  Moreover, $\rho_L$ is pseudo-metric on $\mathcal{R}(M)$ and $\rho_L(g_1,g_2)=0$ if and only if $g_1$ is isometric to $g_2$ [See \cite{MR1835418}, Theorem 7.2.4]. Thus, it can induce a Lipschitz-metric $\rho_{L}$ on $\mathcal{M}(M)$, and then $\rho_{L}$ induces the Lipschitz-topology on $\mathcal{M}(M)$ called L-topology. Furthermore, Lipschitz convergence implies Gromov-Hausdorff Convergence, where the convergence means cauchy sequence convergence relative to their metrics [See \cite{MR2307192}, Proposition 3.6]. 
\begin{proposition}
 L-topology is finer than GH-topology on $\mathcal{M}(M)$.
\end{proposition}

The GH-topology and L-topology on $\mathcal{M}(M)$  only catch the metric information of the basic manifold and  loses much essential information of the smooth manifold. So it may be useful to modify the definition of L-topology on $\mathcal{M}(M)$ to finer topology on $\mathcal{M}(M)$.\par
 For any homorphism of metric space $f:(X,d_{X})\rightarrow(Y,d_{Y})$, the Lipschitz constant of $f$ is defined by $ L(f):= \inf\{ k\geq 1\mid\frac{d_X (x,y)}{k}\leq d_Y (f(x),f(y))\leq kd_X(x,y), x,y\in X \}$. If the set is empty, let $L(f)$ be infinity. 
\begin{lemma}
 Suppose that $M$ and $N$ are smooth closed Riemannian manifolds, then any diffeomorphism of $M$ and $N$ has bounded Lipschitz constant. 
\end{lemma}
\begin{remark}
  The normal of tangent map of diffeomorphism on the unit tangent bundle over closed manifold are uniform bounded, since the tangent map are continuous and the total space of unit tangent bundle over compact manifold are compact.
\end{remark}

 For the composition of diffeomorphism $f\circ g: M \rightarrow N\rightarrow W$,  $L(f\circ g)\leq L(f)\cdot L(g)$ by direct computation. 
 
 \begin{definition}
 Let $g_1, g_2$ in $\mathcal{R}(M)$, $\rho_{SL} (g_1, g_2)=\rho_{SL} ((M,d_1), (M,d_2)):=\inf\{ \log L(f)\mid f\in\mathrm{Diff}\}$, where $Diff$ is the diffeomorphism group of $M$. 
\end{definition}

\begin{lemma}
 $\rho_{SL}$ is a pseudo-metric on $\mathcal{R}(M)$ and  $\rho_{SL}(g_1, g_2)=0$ if and only if $g_1$ is isometric to $g_2$ on $M$.

\end{lemma}
\begin{remark}
 The isometry map between $(M, d_1)$ and $(M, d_2)$ can be constructed by using the closeness of $M$ and the Arzela-Ascoli theorem, if  $\rho_{SL}(d_1, d_2)=0$.
\end{remark}

Continuing the game, one  can define the metric topology on $\mathcal{M}(M)$  called SL-topology by the metric $\rho_{SL}$.
\begin{theorem}
 SL-topology $\preceq$ L-topology $\preceq$ GH-topology $\cong\varepsilon$-topology.
\end{theorem}

Usually those four metrics are not complete metrics on $\mathcal{M}(M)$, so $\mathcal{M}(M)$ is  local compact topology space endowed with their induced metric topology in general. But if we restrict it to the subset of $\mathcal{M}(M)$, it may have some precompact propositions. For example, Gromov precompactness theorem and other convergence theorem on the moduli space [\cite{MR2307192}, Chapter 5].\par

For the non-compact case,  one can ask what is the right topology on $\mathcal{R}^{\geq 0}(V)$ and $\mathcal{M}^{\geq0}(V)$, where $V$ is a non-compact manifold, $\mathcal{R}^{\geq 0}(V)$ is the Riemannian metric with non-negative sectional curvature,  $\mathcal{M}^{\geq 0}(V)$ is the moduli space of $V$ with non-negative sectional curvature?

$\mathbf{Acknowledgment}$: I thank Xuchao Yao for useful discussions. 

\bibliographystyle{alpha}
\bibliography{reference}
\Addresses

\end{document}